\documentclass[12pt]{article}
\usepackage[usenames,dvipsnames,svgnames,table]{xcolor} 
\usepackage{kotex}
\usepackage{graphicx}
\usepackage{epsfig}
\usepackage{amssymb, amsmath, mathtools}
\usepackage{verbatim}
\usepackage{natbib}
\usepackage[affil-it]{authblk}
\usepackage{mathrsfs}
\usepackage{here}
\usepackage{makecell}
\usepackage{tikz}
\usepackage[shortlabels]{enumitem}
\usepackage{yfonts}
\usepackage{comment}
\usepackage{longtable}
\usepackage{tikz,amsmath, amssymb,bm,color}
\usetikzlibrary{circuits.logic.US,circuits.logic.IEC,fit}
\usepackage{algorithm2e}
\usepackage{multirow}
\usepackage{placeins}
\usepackage[colorlinks]{hyperref} 
\usepackage[colorinlistoftodos, textsize=scriptsize]{todonotes} 

\usepackage[framemethod=TikZ]{mdframed}
\mdfdefinestyle{JL}{%
    linecolor=blue,
    outerlinewidth=1pt,
    roundcorner=2pt,
    innertopmargin=0.1\baselineskip,
    innerbottommargin=0.1\baselineskip,
    innerrightmargin=2pt,
    innerleftmargin=2pt,
    backgroundcolor=orange!100!white}

\hypersetup{
  colorlinks,
  citecolor=Black,
  linkcolor=Red,
  urlcolor=Blue}

\includecomment{note} 

\newtheorem{theorem}{Theorem}[section]
\newtheorem{lemma}[theorem]{Lemma}
\newtheorem{definition}[theorem]{Definition}
\newtheorem{proposition}[theorem]{Proposition}
\newtheorem{corollary}[theorem]{Corollary}

\newtheorem{remark}[theorem]{Remark}

\DeclareMathOperator{\argmin}{argmin}
\DeclareMathOperator*{\argmax}{argmax}

\everymath{\displaystyle}

\newcommand{\lra}{\longrightarrow}

\newcommand{\se}[1]{\section{#1}}
\newcommand{\sse}[1]{\subsection{#1}}

\newcommand{\be}{\begin{equation}}
\newcommand{\ee}{\end{equation}}
\newcommand{\bea}{\begin{eqnarray*}}
\newcommand{\eea}{\end{eqnarray*}}
\newcommand{\bean}{\begin{eqnarray}}
\newcommand{\eean}{\end{eqnarray}}
\newcommand{\ben}{\begin{enumerate}}
\newcommand{\een}{\end{enumerate}}
\newcommand{\bi}{\begin{itemize}}
\newcommand{\ei}{\end{itemize}}
\newcommand{\brem}{\begin{remark}}
\newcommand{\erem}{\end{remark}}
\newcommand{\bcen}{\begin{center}}
\newcommand{\ecen}{\end{center}}
\newcommand{\bsv}{\begin{semiverbatim}}
\newcommand{\esv}{\end{semiverbatim}}

\newcommand{\bt}{\begin{theorem}}
\newcommand{\et}{\end{theorem}}
\newcommand{\bl}{\begin{lemma}}
\newcommand{\el}{\end{lemma}}
\newcommand{\bd}{\begin{definition}}
\newcommand{\ed}{\end{definition}}
\newcommand{\bc}{\begin{corollary}}
\newcommand{\ec}{\end{corollary}}
\newcommand{\bp}{\begin{proposition}}
\newcommand{\ep}{\end{proposition}}

\newcommand{\bbR}{ \mathbb{R}}
\newcommand{\bbX}{ \mathbb{X}}

\newcommand{\calB}{\mathcal{B}}
\newcommand{\calC}{\mathcal{C}}

\newcommand{\calF}{\mathcal{F}}

\newcommand{\calL}{\mathcal{L}}
\newcommand{\calM}{\mathcal{M}}

\newcommand{\calR}{\mathcal{R}}

\title{ Post-Processed Posteriors for Banded Covariances}
\author[1]{Kwangmin Lee}
\author[2]{Kyoungjae Lee}
\author[1]{Jaeyong Lee}
\affil[1]{Department of Statistics, Seoul National University}
\affil[2]{Department Statistics, Inha University}

\begin{document}

\maketitle

\begin{abstract}
	We consider Bayesian inference of banded covariance matrices and propose a post-processed posterior.
	The post-processing of the posterior consists of two steps. In the first step, posterior samples are obtained from the conjugate inverse-Wishart posterior which does not satisfy any structural restrictions. In the second step, the posterior samples are transformed to satisfy the structural restriction through a post-processing function.
	The conceptually straightforward procedure of the post-processed posterior makes its computation efficient and can render interval estimators of functionals of covariance matrices. We show that it has nearly optimal minimax rates for banded covariances among all possible pairs of priors and post-processing functions. 
	Furthermore, we prove that, the expected coverage probability of the $(1-\alpha)100\%$ highest posterior density region of the post-processed posterior is asymptotically $1-\alpha$ with respect to a conventional posterior distribution. It implies that the highest posterior density region of the post-processed posterior is, on average, a credible set of a conventional posterior.
	The advantages of the post-processed posterior are demonstrated by a simulation study and a real data analysis.
\end{abstract}

\section{Introduction}
\label{sec:intro}

In this paper, we propose a new Bayesian procedure for banded covariance matrices. The banded matrices are the matrices whose entries farther than a certain distance from the diagonal are all zeros. 
Banded covariance matrices arise  in modelling  marginal dependence structures of variables with natural ordering such as time series data. 
The banded sample covariance  has been applied to the autoregressive and moving average  models  \citep{wu2009banding} and the time-varying autoregressive-moving-average models \citep{wiesel2012covariance}.

When $p$ is small relative to $n$, the inverse-Wishart prior is the most commonly used conjugate prior for the covariance of the multivariate normal model. 
We denote $\Sigma\sim IW_p(\Lambda,\nu)$, if it has density $\pi(\Sigma)\propto |\Sigma|^{-\nu/2} \exp\{-  tr(\Sigma^{-1}\Lambda)/2\}$, for any $p\times p$ positive definite matrix $\Sigma$,   where $\nu >2p$ is the degree of freedom, and  $|\Sigma|$ is the determinant of $\Sigma$.  
The inverse-Wishart prior has many nice properties under the traditional setting of a small $p$. The posterior induced by the inverse-Wishart prior attains the optimal minimax rate when $p \leq c n, ~ 0 \leq c < 1$, under the spectral norm \citep{lee2018optimal}.
The Jeffreys prior for covariance matrices \citep{yang1996catalog} can be expressed as the limit of the inverse-Wishart prior as the degree of freedom and the scale matrix converge to $p+1$ and the $p\times p$ zero matrix, respectively. When the degree of freedom is $2p+1$ and the scale matrix is a diagonal matrix, the marginal distribution of each correlation induced by the inverse-Wishart prior follows a uniform distribution over the interval $[-1,1]$  \citep{huang2013simple}; thus it can be viewed as a non-informative prior for correlations.


When $p \geq n$, however, \cite{lee2018optimal} showed that the degenerate prior $\delta_{I_p}$, an obviously inadequate prior, attains the optimal minimax rate, implying that without further assumptions, the inference of the covariance matrix is hopeless.
This is  expected   because, without any constraint, the number of parameters, $p(p+1)/2$, in the covariance matrix  is   much larger than the sample size $n$.
To reduce the number of effective parameters, several  matrix classes have been proposed including bandable matrices \citep{cai2010optimal,banerjee2014posterior}, sparse matrices \citep{cai2011adaptive,cai2016estimating,lee2019minimax} and  low-dimensional structural matrices \citep{cai2013sparse,pati2014posterior,gao2015rate}.

In this paper, we focus on the banded covariance assumption. The banded covariance assumption is a popular structural assumption to reduce the number of effective parameters, especially when there is a natural ordering between variables. From the frequentist side, the banded covariance estimator has been studied extensively. 
Since the banded covariance structure is an example of the Gaussian covariance graph model, the methods of \cite{kauermann1996dualization} and \cite{chaudhuri2007estimation} can be used for the estimation of the banded covariance. 
However, the two methods are originally designed for  the case of $p < n$, and need a modification of the sample covariance matrix for the case of $p \geq n$. 
\cite{bickel2008regularized} focused on the bandable covariance structure and obtained the convergence rate of the banded sample covariance. 
Despite of a few point estimation methods, there is no frequentist interval estimation method in high-dimensional settings. 
\cite{chaudhuri2007estimation} suggested an interval estimation of  banded covariances under the asymptotic normality of the maximum likelihood estimator, which is valid only for fixed $p$. 
Also, there is no minimax rate result in the literature for banded covariance matrices, although \cite{cai2010optimal, cai2012minimax} showed the tapering estimator satisfies the optimal minimax rate for bandable covariances, which are the matrices whose entries are getting smaller as they are more distant from the diagonal.

Compared  to  frequentist methods, Bayesian methods have a natural advantage of producing interval estimators automatically. 
However, Bayesian methods for banded covariance matrices that are scalable and supported by theoretical properties in high-dimensional settings are scarce.  This is due to the difficulty of  inventing a tractable prior distribution on the space of banded covariances. 
\cite{khare2011wishart} and \cite{silva2009hidden} proposed prior distributions for covariance graphical model, which can be used for banded covariance matrices, but there are no minimax optimality results for these methods.  This is partly because there are no closed forms of normalizing constants for these priors, which prevents direct investigation of posterior asymptotics. 
It is also  mathematically challenging to apply traditional posterior consistency and contraction theorems \citep{ghosal2017fundamentals}, which are applicable when it is hard to tract posterior directly.

In summary, there are no  Bayesian or frequentist methods for banded covariance matrices, which 
(1) are computationally efficient, (2) produce interval estimators for functionals of covariance matrices, and (3) have optimal or nearly optimal  minimax rate.
In this paper, we propose a new Bayesian method that has the above three properties. In particular, we propose
{\it post-processed posteriors} for banded covariance matrices.

The construction of the post-processed posterior consists of two steps, {\it the initial posterior computing step} and {\it the post-processing step}. In the initial posterior computing step, posterior samples are generated from the initial posterior, the conjugate inverse-Wishart posterior for covariance, without any structural restrictions. In the post-processing step, the initial posterior samples are transformed through a function $f(\Sigma)$ whose range belongs to a space of banded covariances. We call the distribution of the transformed posterior samples the post-processed posterior, which will be rigorously defined in Section  \ref{sec:banded_PPP}.

The idea of transforming the posterior samples has been suggested in various settings. Posterior projection methods (Patra and Dunson 2018; Dunson and Neelon 2003; Gunn and Dunson 2005; Lin and Dunson 2014)\nocite{dunson2003bayesian,gunn2005transformation,lin2014bayesian,patra2018constrained} are proposed for various problems, which project the posterior samples onto the  the constrained parameter space to obtain the projected posterior.
Our proposal is the same as the posterior projection method in spirit, but the choice of posterior transformation is determined through asymptotic consideration, while the posterior projection method uses the projection on the constrained space. In fact, our proposal is the posterior projection method on the space of banded covariances with the Frobenius norm.
Recently, \cite{bashir2018post} proposed a support recovery method for sparse precision matrices based on post-processing of the posterior samples.

The post-processed posterior is conceptually straightforward and  computationally  fast. This is advantageous when the data set is huge and the dimension of the observations is high. The existing Bayesian method can be slow at times especially in high-dimensional settings.
Through the simulation study, we will show that the post-processed posterior significantly reduces the computation time compared to the covariance graphical models proposed by \cite{silva2009hidden} and \cite{khare2011wishart}.
Furthermore, the post-processed posterior attains the nearly optimal minimax rate for the class of banded covariance matrices. 
This is the first minimax result for banded covariance matrices in both Bayesian and frequentist sides. 

The banded covariance matrices have been investigated as a case of covariance graphical model, but the minimax lower bound for covariance graphical model is absent in the literature. Methods for obtaining minimax lower bound, e.g.,  Le Cam's method and Assouad's lemma, are  based on the testing problem of $\delta$-separated sets as described in Proposition 15.1 of \cite{wainwright2019high}. Since patterns of parameter spaces of covariance graphical model differ by the graphs, it is not easy to choose representative separated sets for arbitrary graphical structures. 
Instead, we focus on the banded covariance structure and could choose appropriate separated sets.
We also show that the post-processed posterior has the nearly optimal minimax rate for the class of bandable covariances,  which is given in the supplementary material.

It is worth mentioning that there are substantial differences between banded covariance and precision matrices. 
For banded precision matrices, $G$-Wishart priors \citep{banerjee2014posterior} or banded Cholesky priors \citep{lee2017estimating} can be used, and the normalizing constants are available in a closed form. 
Intuitively, in the Bayesian framework, constraints on precision matrices are more manageable than those on covariance matrices because the precision matrix is a natural parameter of multivariate normal distributions as an exponential family. In other words,  the likelihood function of the covariance is expressed through the precision matrix. 
Thus, Bayesian banded covariance matrix estimation is  more challenging than banded precision matrix estimation.  

There is difference in the estimation methods of  covariance and precision matrices in the frequentist literature as well.  The sparse covariance estimation is typically based on banding or thresholding the sample covariance  \citep{cai2012optimal} while the sparse precision matrix estimation is often based on the penalized likelihood approach (Cai et al. 2011; Zhang and Zou 2014)\nocite{zhang2014sparse,cai2011constrained}. 
The difference is due to the form of  the likelihood function as well as the singularity of  the sample covariance matrix. Contrary to the sample covariance matrix, the sample precision matrix,  the inverse of the sample covariance matrix, does not exist when $p>n$, which prevents thresholding the sample precision matrix. One could choose a small constant $\epsilon>0$ to make $S_n + \epsilon I_p$ invertible and use $(S_n + \epsilon I_p)^{-1}$  instead of $S_n^{-1}$ as the sample precision matrix, but it can be computationally unstable especially when $\epsilon$ is  small.



The rest of the paper is organized as follows.
In Section \ref{sec:banded_PPP},  the post-processed posterior is introduced for the banded covariances. 
In Section \ref{sec:minimax}, it is shown that the banding post-processed posterior attains the nearly optimal minimax rate for banded covariance matrices, and  the expected coverage probability of the $(1-\alpha)100\%$ highest posterior density region of the post-processed posterior is asymptotically $1-\alpha$ with respect to a conventional posterior distribution. In Section \ref{sec:numerical}, the post-processed posterior is demonstrated via simulation studies and a real data analysis. 
The supplementary material contains the proofs of the theorems in the paper, a minimax result of for bandable covariance matrices  and more numerical studies.

\section{Post-Processed Posterior}\label{sec:banded_PPP}

Suppose $X_1,  \ldots, X_n$ are independent and identically distributed samples from $N_p(0, \Sigma)$, the $p$-dimensional normal distribution with zero mean vector and covariance matrix $\Sigma =(\sigma_{ij})>0$. 
We write $B >0$ $(B \ge 0)$ if $B$ is a positive (nonnegative) definite matrix.  
When the variables have a natural ordering such as time or causal relationship, it is commonly assumed that the covariance satisfies a band structure.
In this paper, we assume that $\Sigma$ is banded:
\bean
\Sigma \in \mathcal{B}_{p,k} &:=& \mathcal{B}_{p,k}(M_0,M_1) \nonumber \\
&=& \Big\{\Sigma \in \calC_p: \sigma_{ij}=0 \text{ if }|i-j|>k,\forall i, j \in [p],\lambda_{\max}(\Sigma)\le M_0,\lambda_{\min}(\Sigma)\ge M_1 \Big\}, \quad\,\,
\eean
where $0 < M_1 \leq M_0 < \infty$, $k$ is a natural number, $[p] = \{1,2, \ldots, p  \}$, $\calC_p$ is the set of all $p\times p$ positive definite matrices, and $\lambda_{\min}(\Sigma)$ and $\lambda_{\max}(\Sigma)$ are the minimum and maximum eigenvalues of $\Sigma$, respectively.

We propose a computationally efficient and theoretically supported Bayesian method for banded covariance matrices.
The proposed method consists of two steps: the initial posterior computing step and the post-processing step. 
We describe these two steps in detail below.

\begin{enumerate}
	\item[Step 1.] (Initial posterior computing step) \\
	In the initial posterior step, a conjugate posterior for the parameter space without any structural restriction is obtained.
	We take the inverse-Wishart prior $ IW_p(B_0,\nu_0)$.
	We say this is the {\it initial prior} $\pi^i$ for $\Sigma$.
	By conjugacy, the {\it initial posterior} is then
	$$\Sigma  \mathbb{X}_n  \sim IW_p(B_0 + n S_n,\nu_0 +n) ,  $$
	where $\bbX_n =(X_1,\ldots, X_n)^T$ and  $S_n = n^{-1} \sum_{i=1}^n X_i X_i^T $ is the sample covariance matrix.
	We sample $\Sigma^{(1)}, \Sigma^{(2)}, \ldots, \Sigma^{(N)}$ from the initial posteior, $\pi^i(\Sigma | \bbX_n)$.
	\item[Step 2.] (Post-processing step) \\
	Let the function $B_k(B)$ denote the $k$-band operation,
	$$B_k(B) =  \{ b_{ij} I(|i-j| \leq k) \} $$
	for any $B =(b_{ij})\in \bbR^{p\times p}$.
	In the second step, we post-process the samples from the initial posterior to obtain those from the post-processed posterior.
	The samples from the post-processed posterior, $\Sigma_{(i)}$'s, are defined by
	\bean\label{PPP}
	\Sigma_{(i)} &=& f(\Sigma^{(i)}) = B_k^{(\epsilon_n)}(\Sigma^{(i)}) \\
	& := &
	\begin{cases}
		B_k(\Sigma^{(i)}) +  \Big[ \epsilon_n - \lambda_{\min}\{B_k(\Sigma^{(i)})\} \Big] I_p, &\quad\text{ if } \lambda_{\min}\{B_k(\Sigma^{(i)})\}<\epsilon_n, \\
		B_k(\Sigma^{(i)}),   &\quad\text{ otherwise},\nonumber
	\end{cases}
	\eean
	where  $\epsilon_n$ is a small positive number decreasing to $0$ as $n \to \infty$, for $i=1,\ldots,N$. There is no guarantee that  $B_k(\Sigma^{(i)})$ is positive definite, so the second term of \eqref{PPP} is added to make $\Sigma_{(i)}$ positive definite.
	The resulting post-processed samples, $(\Sigma_{(1)},\ldots, \Sigma_{(N)})$, are banded positive definite matrices.
	We suggest using the samples from the post-processed posterior for Bayesian inference of banded covariance matrices.
\end{enumerate}

We call the posterior distribution of \eqref{PPP} the $k$-banding post-processed posterior to emphasize that the $k$-band operation $B_k$ is used; however, other operations can be used to obtain the desired structure.
We call the function $f$ represented by \eqref{PPP} the {\it post-processing function}, and the post-processed posterior with the post-processing function $f$ is denoted by $\pi^{pp}(\cdot | \bbX_n; f)$ or simply $\pi^{pp}(\cdot |\bbX_n)$ if $f$ is understood in the context.




\section{Properties of Post-Processed Posterior}\label{sec:minimax}
\subsection{Minimax Convergence Rates}

In this section, we show that the proposed post-processed posterior procedure is nearly optimal in the minimax sense among all possible post-processed posterior procedures, the pairs of initial priors and post-processing functions.
A conventional Bayesian procedure can be considered as a post-processed posterior procedure, one with prior with support on $\mathcal{B}_{p,k}$ and identity post-processing function. Thus, the proposed post-processed posterior is nearly optimal even compared with conventional Bayesian procedures.

\cite{lee2018optimal} proposed a decision-theoretic framework for comparison of priors. In this framework, a posterior and the space of all probability measures on the parameter space are considered as an action and the action space, respectively. A prior is a decision rule in this setting because a prior combined with data generates a posterior. The posterior-loss (P-loss) and posterior-risk (P-risk) \citep{lee2018optimal} are the loss and risk functions.

The decision-theoretic framework of \cite{lee2018optimal} can be modified for the study of the minimax properties of post-processed posterior. In this setting, a post-processed posterior is an action and a post-processed posterior procedure, a pair of an initial prior and a post-processing function, is a decision rule. We define
the P-loss and P-risk of the post-processed posterior as follows:
\bea
\calL\{\Sigma_0, \pi^{pp}(\cdot \mid \bbX_n; f) \}  & := & E^{\pi^{pp}} ( ||\Sigma_0 - \Sigma || \mid   \bbX_n) \\
& = &  E^{\pi^i} \{ ||\Sigma_0 - f(\Sigma) || \mid \bbX_n\}  ,  \\
\calR(\Sigma_0, \pi^{pp}) & := & E_{\Sigma_0} [\calL\{\Sigma_0, \pi^{pp}(\cdot \mid \bbX_n; f) \} ] \\
& = & E_{\Sigma_0} [ E^{\pi^i}  \{ ||\Sigma_0 - f(\Sigma) || \mid \bbX_n\}  ], 
\eea
where $E^{\pi^i} $ and $E_{\Sigma_0}$ denote  expectations with respect to $\Sigma \sim \pi^i$ and random samples $X_1,\ldots , X_n$ from $N_p(0, \Sigma_0)$, respectively, $\Sigma_0$ is the true value of the $\Sigma$, and $||A|| :=  \{\lambda_{max}(AA^T)\}^{1/2}$ is the spectral norm of a symmetric matrix $A$. 
We now define the minimax rate and convergence rate for post-processed posteriors.
Let
$$\Pi^* = \{ \pi^{pp}(\cdot; f) = (\pi, f): \pi \in \Pi, f \in \calF \}$$
be the space of all possible post-processing procedures, where $\Pi$ is the space of all priors on $\calC_p$, and $\calF$ is the space of all possible post-processing functions, for example, $\calF_k^B  = \{f : \calC_p \to  \calB_{p,k} \}$.

Before we give some definitions of minimax rates, we introduce some notation.
For any positive sequences $a_n$ and $b_n$, we denote $a_n = o(b_n)$ if $a_n / b_n \lra 0$ as $n\to\infty$, and $a_n \lesssim b_n$ if there exists a constant $C>0$ such that $a_n \le C b_n$ for all sufficiently large $n$.
We denote $a_n \asymp b_n$ if $a_n \lesssim b_n$ and $b_n \lesssim a_n$.

A sequence $r_n$ is said to be the minimax rate for $\Pi^*$ over $\calB_{p,k}$ if
$$\inf_{ (\pi, f) \in \Pi^*} \sup_{\Sigma_0 \in  \calB_{p,k}} E_{\Sigma_0}[ \calL\{\Sigma_0, \pi^{pp}(\cdot | \bbX_n; f) \}] \asymp r_n,$$
and a post-processing procedure $(\pi, f) \in \Pi^*$ is said to have P-risk convergence rate $a_n$ if
$$\sup_{\Sigma_0 \in  \calB_{p,k}} E_{\Sigma_0}[ \calL\{\Sigma_0, \pi^{pp}(\cdot | \bbX_n; f) \}] \lesssim a_n.$$
If $a_n \asymp r_n$ and $r_n$ is the P-risk minimax rate, $(\pi, f) \in \Pi^*$ is said to attain the P-risk minimax rate.

We are now ready to state that the banding post-processed posterior attains nearly minimax rate in terms of the P-risk over banded covariance matrices.
Suppose that we observe the data $X_1,\ldots,X_n$ from $p$-dimensional normal distribution, $N_p(0,\Sigma_0)$ with $\Sigma_0 \in \calB_{p,k}$.
The following theorems say that the P-risk of the banding post-processed posterior is nearly minimax optimal.

\begin{theorem}\label{bandupperspec2}
	Let the prior $\pi^i$ of $\Sigma$ be $IW_p(A_n,\nu_n)$. 
	If $A_n\in\calB_{p,k}$ and $n/4\ge (M_0^{1/2}M_1^{-1}\log p) \vee k \vee ||A_n|| \vee (\nu_n-2p)$, then
	\bea
	\sup_{\Sigma_0\in \mathcal{B}_{p,k}} E_{\Sigma_0} \{E^{\pi^i} (|| B_k^{(\epsilon_n)} (\Sigma)-\Sigma_0||^2\mid \bbX_n )\}\le C(\log k)^2\frac{k+\log p}{n},
	\eea
	where the post-processing function $B_k^{(\epsilon_n)}$ is defined in \eqref{PPP}, $\epsilon_n^2 = O\{(\log k)^2 (k+\log p)/n\}$, and $C$ depends on $M_0$ and $M_1$.
	
\end{theorem}

\begin{theorem}\label{theorem:bandlower}
	If $n/2\ge [\min\{(M_0-M_1)^2,1\}\log p] \vee k$, then 
	$$\inf_{ (\pi, f) \in \Pi^*} \sup_{\Sigma_0 \in  \calB_{p,k}} E_{\Sigma_0} \{E^{\pi} (|| f (\Sigma)-\Sigma_0||^2\mid \bbX_n ) \}\ge C\frac{k+\log p}{n},$$
	where $C$ depends on $M_0$ and $M_1$.
	
\end{theorem}

Theorem \ref{bandupperspec2} gives the convergence rate of the P-risk of the banding post-processed posterior for a class of banded covariance matrices $\mathcal{B}_{p,k}$.
A minimax lower bound is given in Theorem \ref{theorem:bandlower}.
The banding post-processed posterior is nearly optimal since its convergence rate has only $(\log k)^2$ factor up to a minimax lower bound.

%
%
%

\subsection{Interval Estimation}

In this subsection, we show  that the $(1-\alpha)100\%$ highest posterior density region of the post-processed posterior is asymptotically on the average an $(1-\alpha)100\%$ credible set of the conventional posterior.
By the conventional Bayesian method, we mean the Bayesian method imposing a prior distribution on banded covariance matrices directly. 
Thus, the  post-processed posterior provides approximations to the  credible regions of the conventional posterior.

For a given integer $0 <  k\le p$ and $\Sigma \in \calC_p$, let $\theta_1=\theta_1(\Sigma)=(\sigma_{ij},|i-j|\le k)$ and $\theta_2=\theta_2(\Sigma)=(\sigma_{ij},|i-j|>k)$.
Let $\pi^c(\theta_1)$ be a prior for $k$-banded covariance matrices. We use the bracket notation for the  distribution or density of random variables. For examples, the joint distribution of $h(X)$ and $g(Y)$ and conditional distribution of $h(X)$ given $g(Y)$  are denoted by $[h(X), g(Y)]$ and   $[h(X)|g(Y)]$, respectively. Probability that $h(X) \in A$ will be denoted by $[h(X) \in A | g(Y)]$ where $A$ is a set.  Subscripts to the brackets are used to distinguish  different joint distributions of $(X, Y)$.

Define 
\begin{eqnarray*}
	\lbrack\theta_1\mid  \bbX_n \rbrack_{PPP, 0} &=& \int \pi^i(\theta_1,\theta_2\mid \bbX_n) d\theta_2,    \\
	&\propto& \int \pi^i(\theta_1,\theta_2) p  \{\bbX_n \mid \Sigma(\theta_1,\theta_2) \} d\theta_2 \\
	\lbrack\theta_1\mid \bbX_n\rbrack_C &=& \pi^c(\theta_1\mid \bbX_n) \\
	&\propto& \pi^c(\theta_1) p  \{\bbX_n\mid \Sigma(\theta_1,0)\}  ,  
\end{eqnarray*}
where $p(\bbX_n\mid \Sigma)$ is the probability density function of  $\bbX_n$ when $X_i$'s follow $N_p(0, \Sigma)$.  
In the above,  $[\theta_1\mid \bbX_n]_{PPP, 0}$ and $[\theta_1\mid \bbX_n]_{C}$ denote the post-processed posterior with only the $k$-band operation $B_k$ and the posterior  of the conventional Bayesian method, respectively.
Note that, in $[\theta_1\mid \bbX_n]_{PPP, 0}$, we use subscript $0$ to distinguish it from the post-processed posterior defined in \eqref{PPP}, which we will denote as $[\theta_1\mid \bbX_n]_{PPP}$.

Suppose that the true covariance matrix $\Sigma_0$ has the $k$-banded structure.
Let $(\hat{\theta}_1^*,\hat{\theta}_2^*)^T =\argmax_{\theta_1,\theta_2} \log p \{\bbX_n \mid \Sigma(\theta_1,\theta_2)\}$ and $\hat{\theta}_1 = \argmax_{\theta_1} \log p \{ \bbX_n \mid \Sigma(\theta_1,0)\}$ be the maximum likelihood estimators. 
Furthermore, we denote the Fisher-information matrix by  
\bea
\mathcal{I}(\theta_1,\theta_2) &=&-E_{\Sigma(\theta_1,\theta_2)} \Big\{ \Big[\frac{\partial}{\partial \theta } \log p \{\bbX_n |\Sigma(\theta_1,\theta_2)\} \Big]^T\Big[\frac{\partial}{\partial \theta } \log p \{\bbX_n |\Sigma(\theta_1,\theta_2) \} \Big] \Big\},  \\
&=& \begin{pmatrix}
	\mathcal{I}_{11} & \mathcal{I}_{12} \\
	\mathcal{I}_{21} & \mathcal{I}_{22} 
\end{pmatrix} 
\eea 
and $\mathcal{I}_{11\cdot 2}(\theta_1,\theta_2) = \mathcal{I}_{11} -\mathcal{I}_{12} \mathcal{I}_{22}^{-1} \mathcal{I}_{21}$.

For Theorem \ref{theorem:coverage}, we assume 
that the total variation distance version of Bernstein von-Mises theorem holds for  $[\theta_1\mid \bbX_n]_{PPP, 0}$  and $[\theta_1\mid \bbX_n]_{C}$, i.e., 

\noindent {\bf A1.} (Bernstein-von Mises condition) 
\bea
\lim\limits_{n\lra \infty}E_{\Sigma_0}||  [ n^{1/2}(\theta_1(\Sigma)-\hat{\theta}^*_1) \mid \bbX_n]_{PPP,0}- N(0,\mathcal{I}_{11\cdot 2}^{-1}  \{\theta_1(\Sigma_0),0 \})||_{TV} &=& 0,\\
\lim\limits_{n\lra \infty}E_{\Sigma_0}|| [n^{1/2}(\theta_1(\Sigma)-\hat{\theta}_1) \mid \bbX_n]_C- N(0,\mathcal{I}_{11}^{-1}\{\theta_{1}(\Sigma_0),0) \})||_{TV} &=& 0.
\eea 
Also, using a slight abuse of notation, we let $N(0,\mathcal{I}^{-1})$ denote the probability measure of the multivariate normal distribution with zero mean vector and covariance matrix $\mathcal{I}^{-1}$.
For any probability measures $P$ and $Q$ on a $\sigma$-field $\calM$, $||P - Q ||_{TV}$ is defined by
$\sup_{A\in \calM} | P(A)- Q(A)|.$
The total variation distance version of the Bernstein von-Mises theorem is given in \cite{van2000asymptotic} and \cite{ghosal2017fundamentals}.


Furthermore, we assume that the following regularity  conditions   hold.  Let $\stackrel{d}{\lra}$ and $\stackrel{P}{\lra}$ denote the convergence in distribution and in probability, respectively.

\noindent {\bf A2.}  As $n\lra\infty$,  
\begin{eqnarray}\label{eq:regconditions}
n^{1/2}L_n' \{ \theta_1(\Sigma_0),0\} &\stackrel{d}{\lra}& N [0,\mathcal{I} \{ \theta_1(\Sigma_0),0\}],\nonumber\\
(\hat{\theta}_1^*,\hat\theta_2^*) &\stackrel{P}{\lra}& \{\theta_1(\Sigma_0),0\}, \\
\hat{\theta}_1 &\stackrel{P}{\lra}& \theta_1(\Sigma_0),\nonumber 
\end{eqnarray}
and $L''_n \{ \theta_1(\Sigma_0),0\} $ is continuous,
where $L_n(\theta_1,\theta_2) =\log p \{\bbX_n \mid \Sigma(\theta_1,\theta_2)\}$, $L_n' \{\theta_1(\Sigma_0),0\} = \partial L_n\{\theta_1(\Sigma_0),0\}  / \partial(\theta_1,\theta_2) $, $L_n''\{\theta_1(\Sigma_0),0\} = \partial^2 L_n\{\theta_1(\Sigma_0),0\} / \partial(\theta_1,\theta_2)^2 $.


Theorem \ref{theorem:coverage} shows that, under the regularity conditions, the highest posterior density region based on the post-processed posterior is on average a credible region of the conventional Bayesian method for banded covariance matrices.

\begin{theorem}\label{theorem:coverage}  Suppose A1 and A2 hold. 
	If $C_{1-\alpha,n}$ is the highest posterior density regions of $[\theta_1 \mid \bbX_n]_{PPP}$ and $p$ is fixed,  then
	\bea
	\lim\limits_{n\lra \infty}E_{\Sigma_0}\{[\theta_1(\Sigma)\in C_{1-\alpha,n} \mid \bbX_n]_C\} = 1-\alpha.
	\eea
	

\end{theorem}

\section{Numerical Studies}\label{sec:numerical}

\sse{A Simulation study: general aspects}

We compare the post-processed posterior with other methods. 
Among frequentist methods, we investigate the performance of banded sample covariance \citep{bickel2008regularized}, dual maximum likelihood estimator \citep{kauermann1996dualization}, and the maximum likelihood estimator by iterative conditional fitting \citep{chaudhuri2007estimation}. 
We also examine the performance of Bayes estimators based on the $G$-inverse Wishart distribution \citep{silva2009hidden} and Wishart distributions for covariance graph \citep{khare2011wishart}. 
Additionally, we conduct the {\it dual post-processed posterior}, which is a post-processing posterior based on the dual algorithm \citep{kauermann1996dualization} instead of the banding post-processing function $B_k$.  
We obtain a posterior sample of the dual post-processed posterior as follows:

\begin{enumerate}
	\item[] Step 1. (Initial posterior computing step) For $l=1,2,\ldots$, sample $\Sigma^{(l)}$ from the initial posterior, 
	$$\Sigma^{(l)} \mid \bbX_n \sim IW_p(B_0+nS_n,\nu_0+n).$$
	
	\item[] Step 2. (Post-processing step) 
	Obtain $\Sigma_{(l)}^{D}$ as the solution of the simultaneous equations:
	$$ \{(\Sigma_{(l)}^{D})^{-1}\}_{ij}  = \{(\Sigma^{(l)})^{-1}\}_{ij}, $$
	for $|i-j|\le k$ and $(\Sigma_{(l)}^{D})_{ij}=0$ for $|i-j|>k$.  
\end{enumerate} 

\begin{figure}[!tb]
	\begin{center}
		\includegraphics[width=14cm,height=10cm]{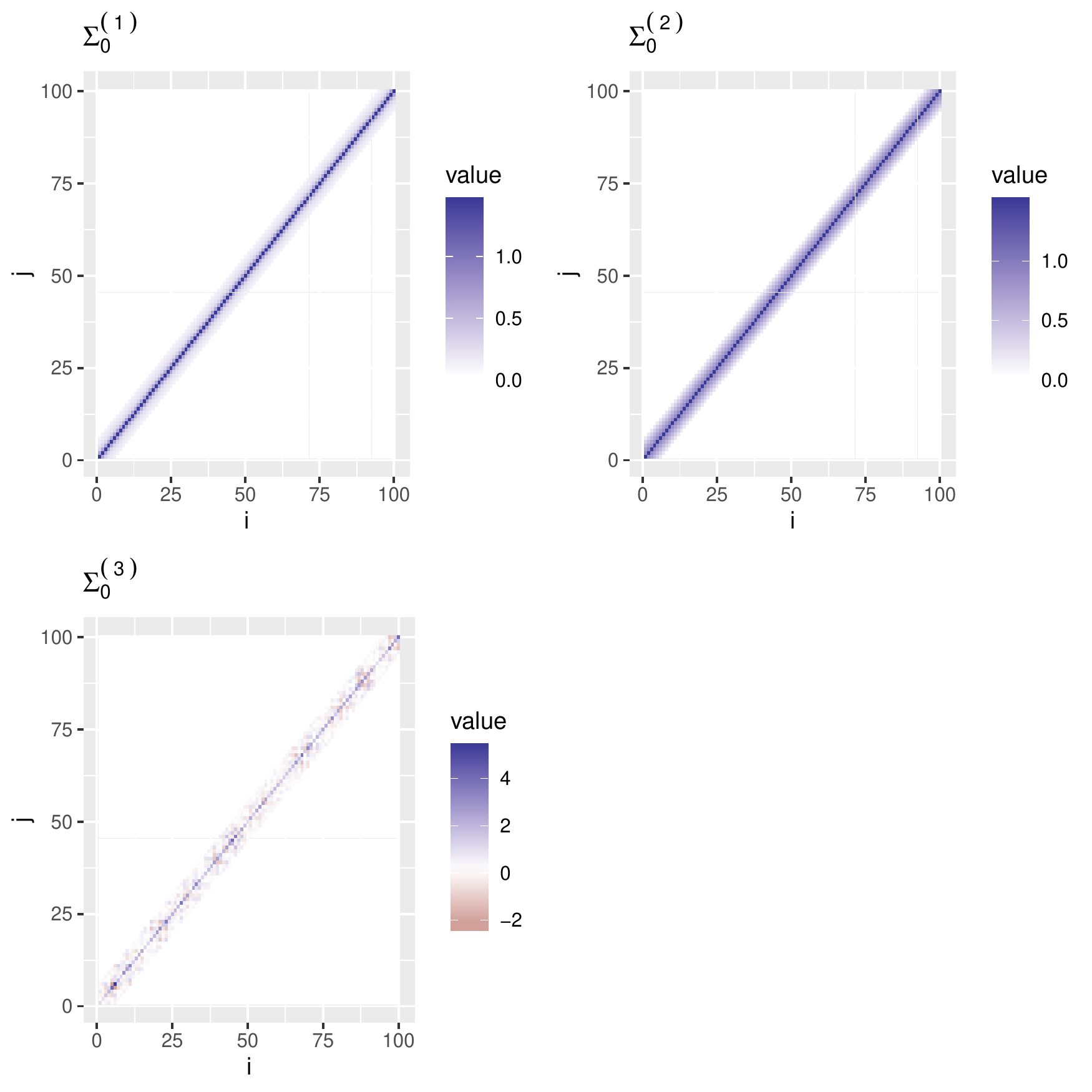}
	\end{center}
	\caption{
		Visualization of true banded covariances.
	}\label{fig:Sigmas} 
\end{figure}

We consider three banded covariances $\Sigma_0^{(1)}$, $\Sigma_0^{(2)}$ and $\Sigma_0^{(3)}$ as the true covariance matrices.
Let  $\Sigma^{(1)*}_0=(\sigma^{(1)}_{0,ij})_{p\times p}$, where
\begin{equation*}
\sigma_{0,ij}^{(1)} = \begin{cases}
1,               & 1\le i=j\le p\\
\rho |i-j|^{-(\alpha+1)}, &  1\le i\neq j\le p,
\end{cases}
\end{equation*}
$\rho=0.6$ and $\alpha=0.1$.  
Then we obtain $\Sigma_0^{(1)}$  by banding $\Sigma_0^{(1)*}$ and adding an identity matrix multiplied by a positive number to make the minimum eigenvalue of resulting matrix to be $0.5$; in particular, $\Sigma_0^{(1)}=B_{k_0}(\Sigma_0^{(1)*}) +[0.5-\{\lambda_{\min}(B_{k_0}(\Sigma_0^{(1)*}))\}]I_p$,  where $k_0$ is the bandwidth. 
Let $\Sigma^{(2)*}_0=(\sigma^{(2)}_{0,ij})_{p\times p}$, where $\sigma_{0,ij}^{(2)} = \{1-|i-j|/(k_0+1) \}\wedge 0$ for any $1 \leq i, j \leq p$. 
Then we set $\Sigma_0^{(2)} = \Sigma^{(2)*}_0 +[0.5-\{\lambda_{\min}(\Sigma_0^{(2)*}) \}]I_p$. Let $\Sigma^{(3)*}_0=L_0D_0L_0^T$ and $\Sigma_0^{(3)} = \Sigma^{(3)*}_0 +[0.5-\{\lambda_{\min}(\Sigma_0^{(3)*}) \}]I_p$, where
\begin{equation*}
L^0_{ij} = \begin{cases}
1,               & 1\le i=j\le p\\
l_{ij}, &   0<i-j \le k_0\\
0, & \text{otherwise},
\end{cases}
\end{equation*}
$l_{ij}$ are independent sample from $N(0,1)$, and $D_0=diag(d_{ii})$ is a diagonal matrix where $d_{ii}$ is independent sample from $IG(5,1)$, the inverse-gamma distribution with the shape parameter $5$ and the scale parameter $1$.  
The true covariance matrices with $p=100$ and $k_0 = 5$ are  plotted   in Figure \ref{fig:Sigmas}.

For each banded covariance with $k_0 = 5$, we generated the data $X_1,\ldots, X_n$ from $N_p(0, \Sigma_0^{(t)})$ independently, where $n=25, 50, 100$ and $p = 100$.
For the  initial prior of the post-processed posterior and dual post-processed posterior, we choose $IW_p(A_0 , \nu_0)$ with $\nu_0 = 2p +3$ and $A_0 = I_p$. 
We use the Bayesian leave-one-out cross-validation method \citep{gelman2014understanding} to choose the adjustment parameter $\epsilon_n$ in the banding post-processing step \eqref{PPP}.   
We define the log-predictive density of $\epsilon_n$ as  
\bea
R(\epsilon_n) &=& \sum_{i=1}^n \log \int p\{X_i \mid B_{k_0}^{\epsilon_n}(\Sigma)\} \pi^i(\Sigma\mid \bbX_{n,-i})d\Sigma \\
&=& \sum_{i=1}^n \log \int p\{X_i \mid B_{k_0}^{\epsilon_n}(\Sigma)\} \frac{\pi^i(\Sigma\mid \bbX_{n,-i})}{\pi^i(\Sigma\mid \bbX_{n})}\pi^i(\Sigma\mid \bbX_{n})d\Sigma,
\eea
where $\bbX_{n,-i}=(X_1,\ldots,X_{i-1},X_{i+1},\ldots, X_n)$ and $p\{\cdot \mid B_{k_0}^{\epsilon_n}(\Sigma)\}$  is the multivariate normal density with zero mean and the covariance matrix $B_{k_0}^{\epsilon_n}(\Sigma)$. Then, using  Monte Carlo method, we obtain the estimated log-predictive density as 
\bean
\hat{R}(\epsilon_n) = \sum_{i=1}^n \log \frac{1}{S}\sum_{s=1}^S p\{X_i \mid B_{k_0}^{\epsilon_n}(\Sigma_s)\} \frac{\pi^i(\Sigma_s\mid \bbX_{n,-i})}{\pi^i(\Sigma_s\mid \bbX_{n})}, \label{formula:BayesCV}
\eean 
and $\hat{\epsilon}_n = \argmin_{\epsilon_n>0}\hat{R}(\epsilon) $,
where $\Sigma_s$ is the $s$th sample of $\pi^i(\cdot \mid \bbX_n)$, and $S$ is the number of the  posterior samples.

For the Wishart distribution for covariance graph \citep{khare2011wishart},  we used $\alpha_i=2k_0+5$ and $U=I_p$ as they suggested.  
Similarly, we set $\delta=5$ and $U=I_p$ for  the $G$-inverse Wishart distribution \citep{silva2009hidden} as they suggested. 
For both methods, the initial values of the $\Sigma$ for the Markov chain Monte Carlo algorithms were set at the identity matrix and 
$500$ posterior samples were drawn.

For the dual maximum likelihood estimator and the maximum likelihood estimator, $S_n + \epsilon_n I_p$ is used  in place of  the sample covariance matrix because these algorithms need a positive definite sample covariance matrix. The adjustment parameter  $\epsilon_n$ is chosen as the minimizer of $\hat{R}_f(\epsilon_n)$, which is defined as 
\bean\label{Rf_epsilon}
\hat{R}_f(\epsilon_n) = \sum_{i=1}^n \log p\{X_i\mid h(\bbX_{n,-i};\epsilon_n)\},
\eean
where  $h(\bbX_n;\epsilon_n)$   is   a frequentist estimator of $\Sigma$ based on $\bbX_n$ and an  adjustment parameter $\epsilon_n$.

For each simulation setting, $100$ sets of samples were  generated.   The performance of each estimator is measured by the  mean spectral norm error 
\bean\label{errordef}
\frac{1}{100} \sum_{s=1}^{100}  ||\Sigma_0-\hat{\Sigma}^{(s)}|| ,
\eean
where $\hat{\Sigma}^{(s)}$ is a point estimate  based on the $s$th simulated data set. For Bayesian methods, we use the posterior mean as the point estimator.
Table \ref{fig:fixederror} shows the   mean spectral norm error   of each method when the true bandwidth $k_0$ is known.  

The maximum likelihood estimator, the banded sample covariance, the Wishart for covariance graph and the post-processed posterior perform well. It appears that the maximum likelihood estimator performs well when  $n/p$ is large, while the post-processed posterior and the banded sample covariance have smaller mean spectral norm errors when  $n/p$ is small.  In all cases, the post-processed posterior  performs reasonably well. 

\begin{table}[!htbp]
	\caption{Specral norm-based errors  for $\Sigma_0^{(1)}$, $\Sigma_0^{(2)}$ and $\Sigma_0^{(3)}$. \label{fig:fixederror}}
	\begin{center}
		\begin{tabular}{lccccccccc}
			
			&\multicolumn{3}{c}{ $n=25$}&\multicolumn{3}{c}{ $n=50$}&\multicolumn{3}{c}{ $n=100$} \\
			& $\Sigma_0^{(1)}$ & $\Sigma_0^{(2)}$ & $\Sigma_0^{(3)}$ & $\Sigma_0^{(1)}$ & $\Sigma_0^{(2)}$ & $\Sigma_0^{(3)}$ & $\Sigma_0^{(1)}$ & $\Sigma_0^{(2)}$ & $\Sigma_0^{(3)}$ \\[5pt]
			Post-processed posterior & 3.67 & 4.62 & 5.63 & 2.16 & 3.01 & 3.61 & 1.48 & 1.94 & 2.34 \\
			G-inverse Wishart & 3.60 & 5.79 & 6.83 & 3.28 & 5.21 & 6.08 & 2.77 & 4.4 & 5.16  \\
			Wishart for covariance graph & 4.56 & 6.85 & 6.08 & 2.07 & 4.36 & 4.81 & 1.41 & 2.9 & 4.96\\
			Dual post-processed posterior & 4.00 & 6.46 & 7.71 & 3.98 & 6.42 & 7.68 & 3.75 & 5.99 & 7.15\\
			Banded sample covariance & 3.38 & 4.5 & 5.66 & 2.19 & 2.8 & 3.42 & 1.51 & 1.9 & 2.23 \\
			Dual maximum likelihood estimator & 3.96 & 6.41 & 7.67 & 3.9 & 6.28 & 7.55 & 3.33 & 5.23 & 6.33 \\
			Maximum likelihood estimator & 4.96 & 4.78 & 6.92 & 2.31 & 2.52 & 3.4 & 1.42 & 1.76 & 2.17
		\end{tabular}
	\end{center}
\end{table}

We compare computation times of the Bayesian methods in Table \ref{fig:computing}. The post-processed posterior is faster than $G$-inverse Wishart distribution and Wishart distribution for covariance graph methods. 
The dual post-processed posterior method is the fastest because it does not have the cross-validation step for the adjustment parameter $\epsilon_n$, but its mean spectral norm errors in Table \ref{fig:fixederror}  shows sometimes poor performance.  

\begin{table}[!htbp]
	\caption{The summary statistics of computing times (unit: sec) for Bayesian methods, when $p=100$ and $n=50$. In the computing times of the post-processed posterior method, the step of Bayesian leave-one-out cross-validation for $\epsilon_n$ is involved.\label{fig:computing}}
	\begin{center}
		\begin{tabular}{lcccccc}
			& $1-$quantile & mean & median & $3-$quantile   \\[5pt]
			Post-processed posterior   & 40.45 & 40.63 & 40.63 & 40.78  \\
			G-inverse Wishart   & 205.47 & 206.67 & 207.32 & 208.23  \\
			Wishart for covariance graph   & 353.91 & 355.14 &  356.31 & 357.08  \\
			Dual post-processed posterior   &  10.60 & 10.73 & 10.67 & 10.78  
		\end{tabular}
	\end{center}
\end{table}

\sse{A Simulation study: interval estimation aspect}

We investigate the performance of interval estimation for functionals of covariances in this section.  
There is no valid frequentist interval estimator for functionals of banded or bandable covariances in the high-dimensional covariance.   But if one assumes $p$ is fixed, the interval estimator for functionals of banded covariances can be derived from the Fisher information matrix given in \cite{chaudhuri2007estimation}. 
Define $vecb(\Sigma):=vecb(\Sigma;k) = vec(\{ \sigma_{ij} : i\leq j , |i-j|\le k\})$ and $Q\in\bbR^{p^2 \times p^*}$ such that $vec(\Sigma)=Q \times vecb(\Sigma;k)$,
where $vec$ is the column-wise vectorization operation, and $p^*$ is the dimension of $vecb(\Sigma;k)$.
By asymptotic normality of maximum likelihood estimators and the Fisher information matrix in \cite{chaudhuri2007estimation}, we obtain
\bea
n^{1/2}\{vecb(\Sigma^{MLE})- vecb(\Sigma_0)\} \stackrel{d}{\lra} N_{p^*}  [0,2 \{Q^T ( \Sigma_0^{-1}\otimes  \Sigma_0^{-1})Q \}^{-1}],
\eea
as $n\lra\infty$, where $\Sigma^{MLE}$ is obtained by the iterative conditional fitting. Let $\phi\{vecb(\Sigma)\}$ and $\nabla\phi\{vecb(\Sigma)\}$ be a functional and its  derivative, respectively. 
By the delta method, we obtain 
\bea
n^{1/2}[\phi\{vecb(\Sigma^{MLE})\}-\phi\{vecb(\Sigma_0)\}]\stackrel{d}{\lra} N(0,\sigma_{0,\phi}^2),
\eea
as $n\lra\infty$, where $\sigma_{0,\phi}^2=2\nabla \phi\{vecb(\Sigma_0)\}\{Q^T ( \Sigma_0^{-1}\otimes  \Sigma_0^{-1})Q \}^{-1}\nabla^T \phi\{vecb(\Sigma_0)\}$. Then, we induce an $(1-\alpha) 100 \%$ confidence interval of the functional as
\bea
\phi(vecb(\Sigma^{MLE})) \pm z_{\alpha/2}\frac{\sigma_{0,\phi}}{n^{1/2}}.
\eea
Since $\sigma_{0,\phi}$ depends on the true covariance matrix, we use an estimated value as
\bea
\hat{\sigma}_\phi = 2\nabla \phi\{vecb(\Sigma^{MLE})\}\{Q^T ( (\Sigma^{MLE})^{-1}\otimes  (\Sigma^{MLE})^{-1})Q \}^{-1}\nabla^T \phi\{vecb(\Sigma^{MLE})\}.
\eea

For Bayesian methods, We obtain credible intervals using the posterior samples. 
For posterior sample $\Sigma_1,\ldots, \Sigma_S$, the $(1-\alpha)100\%$ credible interval for a functional $\phi(\Sigma)$ can be obtained based on $\phi(\Sigma_1),\ldots,\phi(\Sigma_S)$. We set $S=500$ in the simulation.

In the numerical experiment, we focus on the conditional mean for the prediction problem as a functional of covariances.  
When $X_i= (X_{i,1},\ldots,X_{i,p})^T \sim N_p(0,\Sigma)$, the conditional mean given $X_{-p} =(X_1, \ldots, X_{p-1})^T$ is 
\bea
cm(\Sigma; X_{-p}) := E(X_p\mid X_{-p}) = \Sigma_{p,-p} \Sigma^{-1}_{-p,-p} X_{-p}.
\eea
We compare the coverage probabilities and the lengths of intervals for $95\%$ credible intervals of  $cm(\Sigma; X_{-p}) $   in Table \ref{tbl:condmeanprob}.

\FloatBarrier
\begin{table}[!htbp]
	\caption{Coverage probabilities and lengths of interval  estimates  of the conditional mean for banded covariances $\Sigma_0^{(1)}$, $\Sigma_0^{(2)}$ and $\Sigma_0^{(3)}$. The average lengths of intervals are represented in parentheses. \label{tbl:condmeanprob}}
	\begin{center}
		\begin{tabular}{lccccccccc}
			&\multicolumn{3}{c}{ $n=25$}\\
			& $\Sigma_0^{(1)}$ & $\Sigma_0^{(2)}$ & $\Sigma_0^{(3)}$ \\
			Post-processed posterior & 96.7\% (2.54) & 95.5\% (2.27) & 94.3\% (3.66)  \\
			G-inverse Wishart & 44.7\% (1.02) & 49.1\% (0.96) & 45.6\% (1.55)  \\
			Wishart for covariance graph & 99.2\% (2.97) & 99.7\% (3.24) & 97.4\% (3.69)  \\
			Dual post-processed posterior & 75.5\% (0.86) & 62.2\% (0.85) & 46.3\% (1.09)   \\
			Maximum likelihood estimator & 100\% (10.67) & 100\% (13.61) & 100\% (32.88)   \\
			&\multicolumn{3}{c}{ $n=50$}\\
			& $\Sigma_0^{(1)}$ & $\Sigma_0^{(2)}$ & $\Sigma_0^{(3)}$ \\
			Post-processed posterior & 98.3\% (2.05) & 96.7\% (1.93) & 98.5\% (3.32) \\
			G-inverse Wishart & 60.2\% (0.67) & 61.4\% (0.68) & 60.7\% (1.02) \\
			Wishart for covariance graph & 97.7\% (1.49) & 99.3\% (1.85) & 91.5\% (1.8) \\
			Dual post-processed posterior & 80.2\% (0.78) & 72.8\% (0.78) & 58.3\% (1.01) \\
			Maximum likelihood estimator & 100\% (3.12) & 100\% (4.75) & 99.9\% (8.82) \\
			&\multicolumn{3}{c}{ $n=100$}\\
			Post-processed posterior &  95.6\% (1.21) & 96.9\% (1.55) & 98.7\% (2.75)\\
			G-inverse Wishart  & 74.4\% (0.57) & 74.2\% (0.56) & 73.1\% (0.81)\\
			Wishart for covariance graph &  93.3\% (0.92) & 97.3\% (1.11) & 88.3\% (1.09)\\
			Dual post-processed posterior &  50.3\% (0.74) & 53.6\% (0.69) & 49.3\% (1.04)\\
			Maximum likelihood estimator & 99.8\% (1.66) & 100\% (2.84) & 100\% (5.38)  
		\end{tabular}
	\end{center}
\end{table}

The post-processed posterior performs well overall. 
It appears that the post-processed posterior and the Wishart for covariance graph produce practically   reasonable interval estimates.  
When $n = 25$, the post-processed posterior has shorter interval estimates than those of  the Wishart for covariance graph. 
As $n$ increases, the Wishart for covariance graph provides shorter interval estimates, but its coverage probabilities tend to be smaller than the nominal coverage.  
The $G$-inverse Wishart and the dual post-processed posterior have much smaller coverage probabilities than the nominal probability. 
The maximum likelihood estimator tends to produce wide (thus conservative) confidence intervals, which makes it less meaningful in practice.

\sse{Application to call center data}

We apply the post-processed posterior to  analyze  the call center data set, which is used in \cite{huang2006covariance} and \cite{bickel2008regularized}. The data set consists of the number of phone calls for 239 days, and the numbers of calls are recorded for 17 hours from 7:00 and divided into $10$-minute intervals.
We denote the number of calls in the $j$th time index of the $i$th day as $N_{ij}$ ($i=1,\ldots,239; j=1,\ldots,102$), and define $x_{i,j} = (N_{ij} + 1/4)^{1/2}$ so that its distribution is similar to the normal distribution. 
Furthermore, to focus on covariance estimation, we center the data.

Using the covariance estimators by the centered data, we predict the numbers of calls at $j=71,\ldots,102$ time points given those at the other time points. Let $x_i^{(1)} = (x_{i,1},\ldots,x_{i,70})^T$, $x_i^{(2)} = (x_{i,71},\ldots,x_{i,102})^T$, then we obtain estimated conditional mean   of $x_i^{(2)}$ given $x_i^{(1)}$   as
\bea
x_i^{(2)}(\Sigma,x_i^{(1)}) = \Sigma_{21} \Sigma_{11}^{-1} x_i^{(1)} ,
\eea
where $\Sigma_{ab} = E  \{x_{i}^{(a)} (x_{i}^{(b)})^T \}$ for any $a,b \in \{1,2\}$. 
The first 205 days $(i=1,\ldots, 205)$ were used as a training data to estimate $\Sigma$, and the last 34 days $(i=206,\ldots, 239)$ were used as a test data.
We measure accuracy of the methods based on the mean square error,   $ (34)^{-1}\sum_{i=206}^{239} ||x_i^{(2)}-  \hat{x}_i^{(2)}  ||^2$,  where   $\hat{x}_i^{(2)} \equiv x_i^{(2)} (\hat{\Sigma} , x_i^{(1)} )$   is an estimator of ${x}_i^{(2)}$.
Here, $\hat{\Sigma} $ is an estimator of $\Sigma$, where posterior means are used based on 500 posterior samples for Bayesian methods.

Since the true bandwidth is unknown, we choose the bandwidth by the Bayesian leave-one-out cross-validation method , similar to \eqref{formula:BayesCV} but using different $\hat{R}(\cdot)$ as
\bean
\hat{R}(k) : = \sum_{i=1}^n \log \frac{1}{S}\sum_{s=1}^S p\{X_i \mid B_k^{\hat{\epsilon}_n}(\Sigma_s)\} \frac{\pi^i(\Sigma_s\mid \bbX_{n,-i})}{\pi^i(\Sigma_s\mid \bbX_n)} \label{formula:BayesCV2}.
\eean
For the frequentist methods, we select the bandwidth based on the leave-one-out cross-validation similar to \eqref{Rf_epsilon}. 
The mean square error is summarized in Table \ref{tbl:callcenter}.

\begin{table}[!htbp]
	\caption{Mean square error between observations and estimated conditional mean.\label{tbl:callcenter}}
	\begin{center}
		\begin{tabular}{lc}
			Method & error  \\[5pt]
			Post-processed posterior & 0.90 \\
			Inverse-Wishart posterior & 1.22 \\ 
			Dual post-processed posterior & 1.19 \\
			Banded sample covariance & 0.89  \\
			Dual maximum likelihood estimator & 1.03 \\
			Sample covariance & 1.02 
		\end{tabular}
	\end{center}
\end{table}

\begin{figure}[!htbp]
	\begin{center}
		\includegraphics[width=14cm,height=10cm]{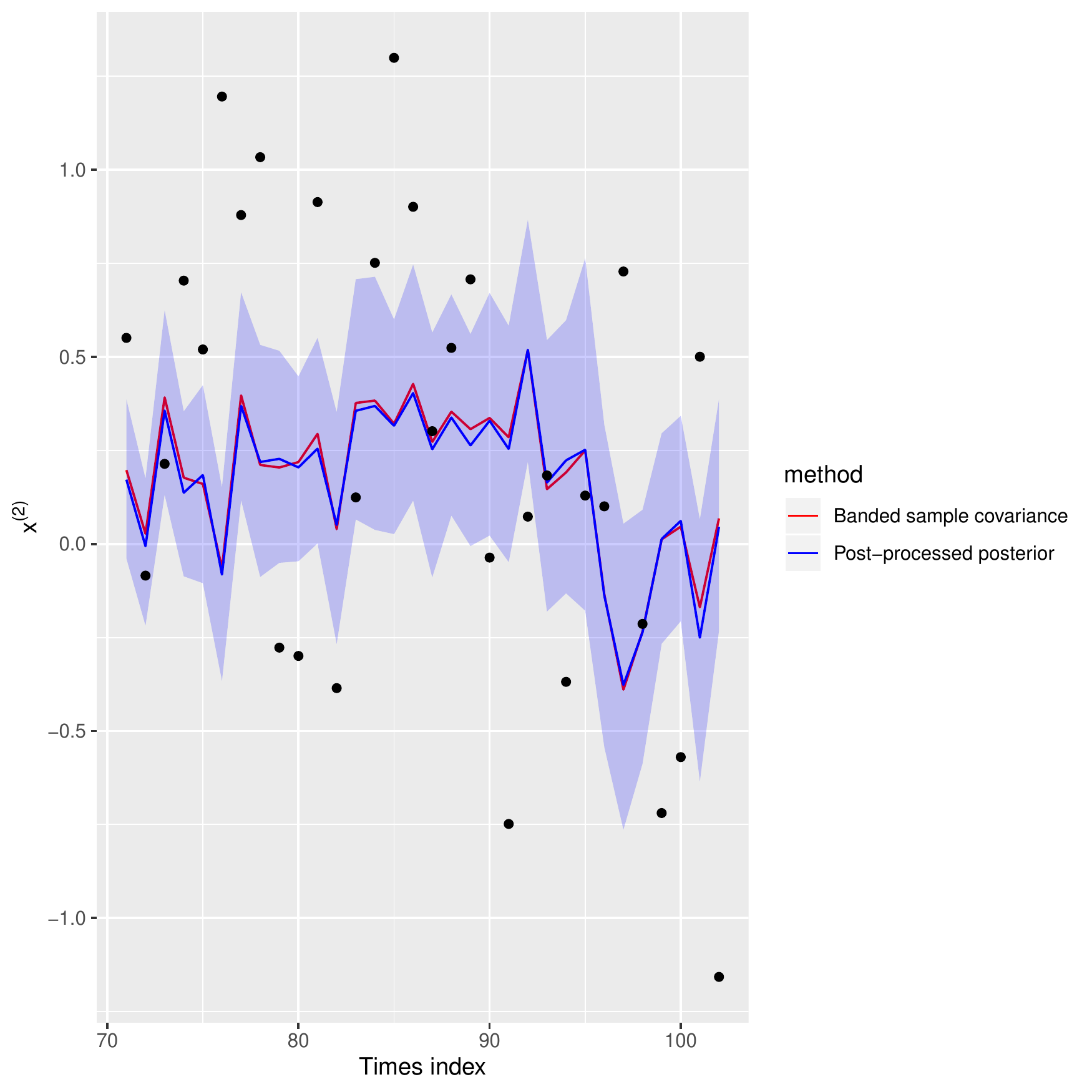}
	\end{center}
	\caption{The estimated conditional mean from the $71$st to $102$nd time indexes of the $2$nd subject. For the post-processed posterior distribution, $95\%$ credible intervals of the conditional mean are also represented as shade.
	}\label{fig:callcenter} 
\end{figure}

The post-processed posterior and banded sample covariance outperform the other methods. 
Although the two methods show similar performance, the post-processed posterior has a benefit over the banded sample covariance because it can give an interval estimator.  
By the definition of $x_i^{(2)}(\Sigma,x_i^{(1)})$, Bayesian methods naturally induce interval estimators based on posterior samples of $\Sigma$.
We visualize the estimators  as well as $95\%$ credible intervals from the post-processed posterior  for the $2$nd subject in the test data in Figure \ref{fig:callcenter}.

\se{Conclusion}
In this paper, we have proposed a non-traditional Bayesian procedure called the post-processed posterior.  It is conceptually straightforward and computationally fast. It attains a nearly minimax convergence rate over all possible pairs of post-processing functions and  initial priors which include  conventional Bayesian posteriors. Also its highest density credible sets are asymptotically  a credible sets of the conventional posteriors on average, and thus its credible sets can be viewed as approximations to the credible sets of the conventional posteriors. 


We applied the post-processing method to the banded covariance structure. But, it can be used in other covariance structures. 
For example, the method can be applied to the class of sparse covariance matrices. We are investigating the theoretical properties of the approach. 
We also believe the post-processing idea can be applied other problems such as sparse linear regression model and high-dimensional nonparametric regression models.  
An open question is to set the boundary of the  post-processing posterior idea: when it has solid theoretical support.

\bibliographystyle{dcu}
\bibliography{cov-ppp}

\end{document}